\documentclass[a4paper,12pt]{article}
\usepackage{amssymb}
\usepackage{amsmath}
\usepackage{epsfig}
\newcommand{\ZZ}{{\mathbb Z}}

\newcommand{\Ga}{\Gamma}
\newcommand{\G}{\Gamma}

\newcommand{\la}{\langle}
\newcommand{\ra}{\rangle}
\newcommand{\qed}{\hfill\hbox{\rule{3pt}{6pt}} \medskip}
\newcommand{\proof}{{\sc Proof. }}

\newcommand{\cay}{\hbox{Cay}}

\newtheorem{theorem}{Theorem}[section]
\newtheorem{lemma}[theorem]{Lemma}

\newtheorem{proposition}[theorem]{Proposition}

\newtheorem{notation}[theorem]{Notation}

\title{ON DISTANCE-REGULAR CAYLEY GRAPHS ON ABELIAN GROUPS}

\author{\v{S}tefko Miklavi\v{c} \\ 
        Andrej Maru\v{s}i\v{c} Institut \\
        University of Primorska \\
        Muzejski trg 2 \\ 
        6000 Koper, Slovenia \\
        stefko.miklavic@upr.si \and
        Primo\v{z} \v{S}parl \\
        Faculty of Education \\
        University of Ljubljana \\
        Kardeljeva plo\v{s}\v{c}ad 16 \\
        1000 Ljubljana, Slovenia \\
        primoz.sparl@pef.uni-lj.si}

\begin{document}
\maketitle

\begin{abstract}
Let $G$ denote a finite abelian group with identity $1$ and let $S$ denote 
an inverse-closed subset of $G \setminus \{1\}$, which generates $G$ and
for which there exists $s \in S$, such that $\la S \setminus \{s,s^{-1}\} \ra \ne G$.
In this paper we obtain the complete classification of distance-regular Cayley graphs 
$\cay(G;S)$ for such pairs of $G$ and $S$.
\end{abstract}

%%%%%%%%%%%%%%%%%%%%%%%%%%%%%%%%%%%%%%%%%%%%%%%%%%%%%%%%%%%%%%%%%%%%%%%%%%%%%%%%%%%%%%%%%%%%
\section{Introduction}
\label{sec:intro}

A connected finite graph is {\em distance-regular} if the cardinality of the 
intersection of two spheres depends only on their radii and the distance between
their centres. Even though this condition is purely 
combinatorial, the notion of distance-regular graphs is closely related to 
certain topics in algebra, and has motivated a development of
various new algebraic notions, as well as shed a new light on the
existing ones. This interplay of concepts proves to be especially intimate
when a subclass of distance-regular {\em Cayley graphs} is considered
(see Section \ref{sec:prelim} for formal definitions of distance-regular graphs and Cayley graphs).

\smallskip
Among distance-regular Cayley graphs, those of diameter $2$ (also called
{\em strongly regular} Cayley graphs) have been investigated most
thoroughly. Such graphs are equivalent to so-called regular partial
difference sets (see \cite{MaSurvey} for the survey of this topic),
and many results on strongly regular Cayley graphs
are formulated in the language of partial difference sets.
Though many authors investigated such graphs, a complete classification
still seems to be beyond reach. In fact, not even strongly regular Cayley graphs of abelian groups have been classified. However, cyclic groups seem to be
easier to handle. Strong regularity of circulants (that is Cayley graphs of cyclic groups) has been investigated
by several authors and a complete classification of strongly regular
circulants was independently achieved by Bridges and Mena \cite{BM}, 
Ma \cite{Ma1}, and partially by Maru\v{s}i\v{c} \cite{DM}.
More recently Muzychuck \cite{Mu} classified strongly regular Cayley graphs on
$\ZZ_{p^n} \times \ZZ_{p^n}$, $p$ a prime.
As for general distance-regular Cayley graphs, 
distance-regular Cayley graphs over cyclic and dihedral groups were classified in \cite{MP, MP1}.

\smallskip 
Let $G$ denote a finite abelian group with identity $1$ and let $S$ denote 
an inverse-closed subset of $G \setminus \{1\}$, which generates $G$ and
for which there exists $s \in S$, such that $\la S \setminus \{s,s^{-1}\} \ra \ne G$. 
The main result of this paper is the following complete classification of distance-regular Cayley graphs $\cay(G;S)$ for such pairs of $G$ and $S$ (see Section~\ref{sec:prelim} for the
definitions of the graphs appearing in the theorem).

\begin{theorem}
\label{the:main}
Let $G$ be an abelian group with identity $1$ and let $S$ be an inverse closed subset of $G \setminus \{1\}$ which generates $G$ and for which there exists 
$s \in S$ such that $\la S \setminus \{s, s^{-1}\} \ra \neq G$. Then the Cayley graph $\cay(G;S)$ is distance-regular if and only if it is isomorphic to one of the following graphs:
\begin{itemize}
\item[(i)] The complete bipartite graph $K_{3,3}$.
\item[(ii)] The complete tripartite graph $K_{2,2,2}$.
\item[(iii)] The complete bipartite graph minus a $1$-factor $K_{6,6} - 6K_2$.
\item[(iv)] The cycle $C_n$ for $n \geq 3$.
\item[(v)] The Hamming graph $H(d,n)$, where $d \geq 1$ and $n \in \{2,3,4\}$.
\item[(vi)] The Doobs graph $D(n,m)$ where $n,m \geq 1$.
\item[(vii)] The antipodal quotient of the Hamming graph $H(d,2)$, where $d \geq 2$.
\end{itemize}
\end{theorem}

Note that the above infinite families include the complete graphs $K_2$, $K_3$ and $K_4$, and that 
the above graphs are distance-regular \cite{BCN}[Section 9.2].

%%%%%%%%%%%%%%%%%%%%%%%%%%%%%%%%%%%%%%%%%%%%%%%%%%%%%%%%%%%%%%%%%%%%%%%%%%%%%%%%%%%%%%%%%%%%
\section{Preliminaries}
\label{sec:prelim}

In this section we review some definitions and basic facts about distance-regular
graphs and Cayley graphs.
More background information on distance-regular graphs
can be found in \cite{BCN}. 

\medskip
Throughout this paper all graphs
are assumed to be finite, undirected and without loops or multiple
edges. For a graph $\Ga$ we let
$V=V(\Ga)$, $E=E(\Ga)$ and $\partial_{\Ga}$ (or just $\partial$)
denote the vertex set, the edge set
and the path length distance function, respectively.
The {\em diameter}
$\max \{\partial(x,y) | x,y \in V(\Ga)\}$ of $\Ga$ will be denoted by $d_{\Ga}$
(or just $d$, when the graph $\Ga$ is clear form the context).
For a positive integer $n$ we denote by $K_n$ the complete graph on $n$ vertices, and by 
$C_n \; (n \ge 3)$ the cycle on $n$ vertices.

\medskip
For a vertex $x \in V(\Ga)$
and an integer $i$ we let $N_i(x) =\{ y\mid \partial(x,y) = i\}$
denote the $i$-th sphere centred at $x$. We abbreviate $N(x)=N_1(x)$.
For a connected graph $\G$ and $x,y \in V$ with $\partial(x,y)=i$ we denote
$$
  c_i(x,y) = |N_{i-1}(x) \cap N(y)|, \qquad 
  a_i(x,y) = |N_i(x) \cap N(y)|, \qquad
  b_i(x,y) = |N_{i+1}(x) \cap N(y)|.
$$
If $c_i(x,y)$ ($a_i(x,y)$, $b_i(x,y)$ respectively) does not depend on the choice of $x,y$ with
$\partial(x,y)=i$ (but only depends on the distance $i$ between $x$ and $y$),
then we say that the {\em intersection number} $c_i$ ($a_i$, $b_i$, respectively) {\em exists} for $\G$ and we set $c_i=c_i(x,y)$ ($a_i=a_i(x,y)$, $b_i=b_i(x,y)$), 
where $x,y \in V$ with $\partial(x,y)=i$.
Observe that if the intersection numbers $a_1, a_2, \ldots , a_i$ all exist, then $a_1 = a_2 = \cdots = a_i = 0$ holds 
if and only if for each $j \in \{3,5, \ldots, 2i+1\}$ there is no cycle 
of length $j$ in $\G$. Observe also that $a_0$, $c_0$, $c_1$ and $b_d$
always exist and $a_0=0$, $c_0=0$, $c_1=1$ and $b_d=0$ holds. 
A connected graph $\Ga$ with diameter $d$ is said to be {\em distance-regular} whenever
the intersection numbers $c_i, a_i$ and $b_i$ exist for all $0 \le i \le d$.
Note that a distance-regular graph $\G$ is regular with valency $k=b_0$, and 
\begin{equation}
\label{eq:int}
  a_i+b_i+c_i=k \qquad (0 \le i \le d).
\end{equation} 
Note also that $b_i \ne 0$ for $0 \le i \le d-1$ and $c_i \ne 0$ for $1 \le i \le d$.
The array
\begin{equation}
\label{intersectionarray}
\{ b_0, b_1, \ldots, b_{d-1}; c_1, c_2, \ldots, c_d\}
\end{equation}
is called the {\em intersection array} of $\Ga$.

Let $\G$ denote a connected graph with diameter $d$. 
For each integer $i \, (0 \le i \le d)$, the $i$th {\it distance matrix} $A_i$
has rows and columns indexed by the elements of $V$, with the $x,y$ entry defined as
\begin{equation}
  (A_i)_{x y} = \left\{ \begin{array}{lll}
                 1 & \hbox{if } \; \partial(x,y)=i,    & \\
                   &                                   & (x,y \in V). \\
                 0 & \hbox{if } \; \partial(x,y) \ne i & \end{array} \right. 
\end{equation}
Note that $A_1$ is the usual adjacency matrix of $\G$. It is easy to see that $\G$ is distance-regular
if and only if there exist nonnegative integers $a_i, b_i$ and $c_i$ such that
$$
  A_1 A_i = b_{i-1}A_{i-1} + a_i A_i + c_{i+1} A_{i+1} \quad  (0 \le i \le d),
$$
where $A_{-1}$ and $A_{d+1}$ are taken to be the zero matrices. Of course, these integers coincide with the intersection numbers from the above definition of distance-regular graphs.

\medskip
For graphs $\Ga_1$ and $\Ga_2$, their {\em Cartesian product} is a graph $\Ga_1 \square \Ga_2$ with vertex-set
$V(\Ga_1) \times V(\Ga_2)$, where vertices $(u_1,v_1), (u_2,v_2) \in V(\Ga_1) \times V(\Ga_2)$ are adjacent 
if and only if $u_1 = u_2$ and $v_1, v_2$ are adjacent in $\Ga_2$, or $u_1, u_2$ are adjacent in $\Ga_1$ and $v_1=v_2$.
It is well known that 
\begin{equation}
\label{eq:dist}
  \partial_{\Ga_1 \square \Ga_2}((u_1,v_1),(u_2,v_2)) = \partial_{\Ga_1}(u_1,u_2) + \partial_{\Ga_2}(v_1,v_2).
\end{equation}
For positive integers $d$ and $q$, the Hamming graph $H(d,q)$ is the Cartesian product of $d$ copies of 
the complete graph $K_q$. Note that in the case of $q = 2$ (these graphs are known as the {\em hypercube graphs}) each vertex has a unique counterpart, the so-called antipodal vertex, at maximal distance $d$. Such a pair of
antipodal vertices is thus a block of imprimitivity for the automorphism group of $H(d,2)$.
The {\em antipodal quotient} of the graph $H(d,2)$ has as vertices the above mentioned blocks consisting of pairs of antipodal vertices with two such blocks adjacent 
whenever there is an edge between these blocks in $H(d,2)$. By~\cite{BCN}[Proposition 4.2.2] the antipodal quotient of $H(d,2)$ is a distance-regular graph.
For a nonnegative integer $n$ and a positive integer $m$, the Doobs graph $D(n,m)$ is the Cartesian product of $H(n,4)$ with 
$m$ copies of the Shrikhande graph~\cite{SS59} (where in the case of $n = 0$ we just take the Cartesian product of $m$ copies of the Shrikhande graph). 
The Doobs graph $D(n,m)$ is distance-regular with the same 
intersection numbers as $H(n+2m,4)$ (see~\cite{BCN}[page 262]).

\medskip 
Let $G$ denote a finite group with identity $1$ and let $S$ denote an inverse-closed subset of
$G \setminus \{1\}$. The {\em Cayley graph} $\cay(G;S)$ of the group $G$ with respect to the 
{\em connection set} $S$ is the graph with vertex-set $G$, 
in which $g \in G$ is adjacent with $h \in G$ if and only if $h=gs$ for some $s \in S$.
Observe that $\cay(G;S)$ is regular with valency $k=|S|$ and is connected if and only if $\la S \ra = G$.

%%%%%%%%%%%%%%%%%%%%%%%%%%%%%%%%%%%%%%%%%%%%%%%%%%%%%%%%%%%%%%%%%%%%%%%%%%%%%%%%%%%%%%%%%%%%
\section{The graphs from the main result}
\label{sec:graphs}

In this section we prove that each of the graphs appearing in Theorem~\ref{the:main} is isomorphic to a Cayley graph $\cay(G;S)$, where
$G$ and $S$ are as in Theorem~\ref{the:main}. However, in this section only, we will be using additive notation for the abelian group $G$.

It is easy to see that  $K_{3,3} \cong \cay(\ZZ_6;\{\pm 1, 3\})$ (and so we can take $s = 1$), $K_{2,2,2} \cong \cay(\ZZ_6 ; \{\pm 1, \pm 2\})$ (and so we can take $s = 1$) 
and $K_{6,6} - 6K_2 \cong \cay(\ZZ_6 \times \ZZ_2 ; \{\pm (1, 0),$ $\pm (2, 1), (0,1)\})$ (and so we can take $s=(1,0)$).
Clearly $C_n \cong \cay(\ZZ_n ; \{\pm 1\})$ (and so we can take $s = 1$). 

It is well known that $H(d,n) \cong \cay((\ZZ_n)^d ; S)$, where $S$ consists of all the elements of the form $(0,0,\ldots,0,x,0,0,\ldots , 0)$, where $x \neq 0$ and the 
position (coordinate) of $x$ is arbitrary. This of course implies that $|S| = d(n-1)$. For $2 \leq n \leq 4$ we can thus take $s = (1,0,0,\ldots , 0)$. Note that in the case of 
$n \in \{2,4\}$ we have that $\la S \setminus \{s, s^{-1}\} \ra$ is an index $2$ subgroup of $G$, while in the case of $n = 3$ this subgroup has index $3$ in $G$.

One can verify that the Shrikhande graph is isomorphic to the Cayley graph $\cay(\ZZ_4 \times \ZZ_4 ; \{\pm (1, 0), \pm (1, 1), \pm (0,1)\})$. (Note that it can be shown that
this graph is not of the form required by Theorem~\ref{the:main}. We only need this description to analyze Doobs graphs $D(n,m)$ for $n, m \geq 1$.)
It is well known and easy to see that for Cayley graphs $\cay(G_1;S_1)$ and $\cay(G_2;S_2)$ of abelian groups $G_1$ and $G_2$ (this holds also for general groups) 
the Cartesian product $\cay(G_1;S_1) \square \cay(G_2;S_2)$ is isomorphic to the Cayley graph $\cay(G_1 \times G_2 ; S)$, where 
$S = \{(x,0), (0,y)\, :\, x \in S_1, y \in S_2\}$. 
Thus every Doobs graph $D(n,m)$ is a Cayley graph of an abelian group. Moreover, since in our case $n \geq 1$, we have that $D(n,m)$ is isomorphic to $D(n-1,m) \square K_4$.
Now, letting $D(n-1,m) \cong \cay(G;S_1)$ we have $D(n,m) \cong \cay(G \times \ZZ_4 ; \{(s_1,0), (0,s_2)\, :\, s_1 \in S_1, s_2 \in \{\pm 1, 2\}\})$ 
(and so we can take $s = (0,1)$).

For $n \geq 2$ the Hamming graph $H(n,2)$ is isomorphic to the Cartesian product $C_4 \square H(n-2,2)$. Let 
$H(n-2,2) \cong \cay((\ZZ_2)^{n-2} ; S_2\})$, where $S_2$ consists of all the elements $e_i = (0,0,\ldots , 0,1,0,0, \ldots , 0)$ with $1$ being on the $i$-th
coordinate. Then $H(n,2) \cong \cay(\ZZ_4 \times (\ZZ_2)^{n-2} ; \{(\pm 1, 0), (0,s_2) \,:\, s_2 \in S_2\})$. The antipodal vertex of the vertex $(x,y)$, 
$x \in \ZZ_4$, $y \in (\ZZ_2)^{n-2}$ is thus $(x+2, y+e)$, where $e = e_1 + e_2 + \cdots + e_{n-2}$. It is now easy to verify that the antipodal quotient of $H(n,2)$ is 
isomorphic to the Cayley graph $\cay(\ZZ_4 \times (\ZZ_2)^{n-3} ; S)$, where $S$ consists of the elements $(\pm 1,0)$, $(2, e')$ and $(0,e_i')$, $1 \leq i \leq n-3$, 
where $e_i' \in (\ZZ_2)^{n-3}$ is the element $(0,0,\ldots , 0,1,0,0, \ldots , 0)$ with $1$ being on the $i$-th coordinate.  Therefore we can take $s = (1,0)$.
\bigskip

The proof that each distance-regular Cayley graph $\G = \cay(G;S)$, where $G$ and $S$ are as in Theorem~\ref{the:main}, is rather complex. 
Various cases need to be considered. 
Clearly, if the valence of $\G$ is less than $3$, we have that $\G \cong K_2$ or $\G \cong C_n$. 
For the rest of the paper we thus only consider graphs of valence at least $3$.
To simplify the statements and to fix the notation used in the arguments, we adopt the following notational convention for the rest of the paper.

\begin{notation}
\label{blank}
Let $G$ denote a finite abelian group with identity $1$ and let $S$ be an inverse-closed subset of
$G \setminus \{1\}$ of cardinality at least $3$. Assume that there exists $s \in S$ such that $\la S \setminus \{s,s^{-1}\}\ra \ne G$
and abbreviate $H = \la S \setminus \{s,s^{-1}\}\ra$. Let $o(s)$ denote the order of the element $s$ and 
let $[G:H]$ denote the index of the subgroup $H$ in $G$. Let $\Ga = \cay(G;S)$ and
let $\Ga_0 = \cay(H; S \setminus \{s,s^{-1}\})$. 
\end{notation}

Observe that $[G:H]$ divides $o(s)$. 
The following result that will be used throughout the paper (at times without explicit reference to it) is straightforward and therefore left to the reader.

\begin{lemma}
\label{lem:easy}
Let $G$ and $\G$ be as in Notation~\ref{blank}. Pick $g,h \in G$.
If $g$ and $h$ are in different cosets of $H$, then they are adjacent 
if and only if $h=gs$ or $h=gs^{-1}$.
\end{lemma}

\begin{lemma}
\label{lem:c2}
Let $\G$ be as in Notation~\ref{blank}. If $\Ga$ is distance-regular, then $2 \leq c_2 \leq 4$ unless 
$\Ga$ is the complete graph $K_4$.
\end{lemma}
\proof
Since $|S| = k \ge 3$, the set $S \setminus \{s,s^{-1}\}$ is not empty. Since $G$ is abelian, it follows that 
for any $t \in S \setminus \{s,s^{-1}\}$ the vertices $s$ and $t$ are common neighbours of $1$ and $st$. By
Lemma~\ref{lem:easy} the vertices $1$ and $st$ are adjacent if and only if $st = s^{-1}$. Thus, 
if $S \setminus \{s, s^{-1}, s^{-2}\}$ is nonempty, we have that $c_2 = |N(1) \cap N(st)| \ge 2$ for any 
$t \in S \setminus \{s, s^{-1}, s^{-2}\}$. Moreover, by Lemma~\ref{lem:easy} we have that $N(1) \cap N(st) \subseteq \{s,s^{-1}, t, s^2 t\}$, and so $c_2 \leq 4$. 
In the only remaining case we have that $S = \{s, s^{-1}, s^{-2}\}$, so that 
$s^2 = s^{-2}$ must hold. It follows that $\Ga \cong K_4$, as claimed.
\qed 

\begin{proposition}
\label{thm:index}
Let $\G$ be as in Notation~\ref{blank}.
If $\Ga$ is distance-regular, then $[G:H] \le 4$. Moreover, if $[G:H] = 4$, then $o(s)=4$.
\end{proposition}
\proof
Suppose $[G:H] \ge 4$ and consider the vertices $1$ and $s^2$ of $\Ga$. Observe that, by Lemma~\ref{lem:easy}, 
$s$ is the unique common neighbour of $1$ and $s^2$ contained in $Hs$. By Lemma~\ref{lem:c2} $c_2 \ge 2$ holds, and so 
there exists at least one more common neighbour of $1$ and $s^2$. By Lemma~\ref{lem:easy} this vertex is on one hand $s^{-1}$
and on the other hand $s^3$. It follows that $s^3 = s^{-1}$, and so $o(s)=4$ and $[G:H]=4$. \qed

\medskip
In what follows we split our analysis into four cases, each dealt with in a separate section:
\begin{itemize}
\item $[G:H] = o(s)$,
\item $[G:H] = 3$ and $o(s) \ge 6$,
\item $[G:H] = 2$ and $o(s) \ge 6$,
\item $[G:H] = 2$ and $o(s) = 4$.
\end{itemize}

%%%%%%%%%%%%%%%%%%%%%%%%%%%%%%%%%%%%%%%%%%%%%%%%%%%%%%%%%%%%%%%%%%%%%%%%%%%%%%%%%%%%%%%%%%%%
\section{Case $[G:H] = o(s)$}
\label{sec:case1}

In this section we consider the distance-regular graphs $\Ga$, where $\Ga$ is as in Notation~\ref{blank}, for which the order of $s$ is equal to the index $[G:H]$. 
The following result is straightforward and therefore left to the reader (recall that $\Ga_0 = \cay(H; S \setminus \{s,s^{-1}\})$).

\begin{lemma}
\label{lem:cartesian}
Let $\G$ be as in Notation~\ref{blank} and suppose that $o(s) = [G:H]$.
Then $\Ga$ is isomorphic to $C_{o(s)} \square \Ga_0$ if $o(s) \ge 3$,
and to $K_2 \square \Ga_0$ if $o(s)=2$.
\end{lemma}

\begin{proposition}
\label{prop:cartesian}
Let $\G$ be as in Notation~\ref{blank} and suppose that $o(s) = [G:H]$ and that 
$\Ga$ is a distance-regular graph with valency $k \ge 3$ and diameter $d$.
Then the following holds.
\begin{itemize}
\item[{\rm (i)}]
If $o(s) \in \{2, 4\}$ then $\G$ is isomorphic to the Hamming graph $H(d,2)$.
\item[{\rm (ii)}]
If $o(s) = 3$ then $\G$ is isomorphic to the Hamming graph $H(d,3)$.
\end{itemize}
\end{proposition}
\proof
(i) Suppose first that $o(s) = [G:H] = 4$. By Lemma~\ref{lem:cartesian}, $\G$ is isomorphic to $C_4 \square \Ga_0$.
Observe that, by \eqref{eq:dist}, the diameter of $\Ga_0$
is $d-2$. We first claim that the intersection numbers $c_i$ of $\Ga$
satisfy $c_i = i$ for $1 \le i \le d$. We show this by induction on $i$. Note that the claim is true for $i=1$.
Since, by Lemma~\ref{lem:easy}, the only common neighbours of $1$ and $s^2$ are $s$ and $s^3$, the claim is true also for $i=2$.
Let now $i \ge 3$. Pick $h \in H$ such that $\partial_{\Ga}(1,h)=i-2$ and note that, by \eqref{eq:dist}, we have 
$\partial_{\Ga}(s^2,h)=i$. Denote $C=N(s^2) \cap N_{i-1}(h)$ and let us compute $c_i = |C|$.
By Lemma~\ref{lem:easy} we have $C = \{s,s^3\} \cup (C \cap Hs^2)$.
However, by Lemma~\ref{lem:cartesian} and \eqref{eq:dist}, $|C \cap Hs^2| = |N(s) \cap Hs \cap N_{i-2}(h)|$.
But by Lemma~\ref{lem:easy} we have $N(s) \cap N_{i-2}(h) = \{1\} \cup (N(s) \cap Hs \cap N_{i-2}(h))$ and therefore 
$|C \cap Hs^2| = c_{i-1} - 1$. It follows that $c_i = c_{i-1} + 1$, which completes the induction step.

We next claim that the intersection numbers $a_i$ of $\Ga$ satisfy $a_i=0$ for $0 \le i \le d$.
Note that $a_0=0$ by definition. By Lemma~\ref{lem:easy}, adjacent vertices $1$ and $s$ have no common
neighbours, and so $a_1=0$. Suppose that $a_{i-1} = 0$ for some $2 \le i \le d$. 
Pick $h \in H$ such that $\partial_{\Ga}(1,h)=i-2$. As above, $\partial_{\Ga}(s^2,h)=i$.
If $a_i \ne 0$, then
there is at least one vertex $g$ in the set $N(s^2) \cap N_i(h)$. By Lemma~\ref{lem:easy} and \eqref{eq:dist}, 
$g$ is contained in the coset $Hs^2$, so $g = h' s^2$ for some $h' \in H$. Consider now the vertex $h' s$ and note that it is in the set
$N(s) \cap N_{i-1}(h)$. But this implies that $a_{i-1} ( = |N(s) \cap N_{i-1}(h)|) \ne 0$, a contradiction.
Therefore, $a_i = 0$ and the claim is proven. 

By \eqref{eq:int} we now find $k=c_d + a_d = d$ and $b_i = k-i = d-i$ for $0 \le i \le d-1$. 
By \cite{BCN}[Corollary 9.2.5], $\Ga$ is the Hamming graph $H(d,2)$.

\smallskip \noindent
Suppose now $o(s) = [G:H] = 2$. By Lemma~\ref{lem:cartesian}, $\G$ is isomorphic to $K_2 \square \Ga_0$.
By \eqref{eq:dist} the diameter of $\Ga_0$ is $d-1$. 
Similarly as above we find that the intersection numbers $a_i$, $b_i$ and $c_i$ of $\Ga$
satisfy $a_i=0$, $b_i=d-i$ and $c_i = i$ for $0 \le i \le d$. 
By \cite{BCN}[Corollary 9.2.5], $\Ga$ is the Hamming graph $H(d,2)$.

\medskip \noindent
(ii) Since $o(s) = [G:H] = 3$, $\G$ is isomorphic to $C_3 \square \Ga_0$.
By \eqref{eq:dist} the diameter of $\Ga_0$ is $d-1$. 
We first claim that the intersection numbers $a_i$ and $c_i$ of $\Ga$
satisfy $a_i = c_i = i$ for $1 \le i \le d$. We show this by induction on $i$. 
Note that $c_1=1$ by definition. By Lemma~\ref{lem:easy}, the only common neighbour of $1$ and $s$ is $s^2$,
so also $a_1 = 1$.
For $2 \le i \le d$ pick $h \in H$ such that $\partial_{\Ga}(1,h) = i-1$ and note that
$\partial_{\Ga}(s,h)=i$. Similarly as in (i) above we find
$$
  c_i = |N(s) \cap N_{i-1}(h)| = |\{1\}| + |N(s) \cap Hs \cap N_{i-1}(h)| = 1 + |N(1) \cap N_{i-2}(h)| = 1 + c_{i-1}
$$
and
$$
  a_i = |N(s) \cap N_i(h)| = |\{s^2\}| + |N(s) \cap Hs \cap N_i(h)| = 1 + |N(1) \cap N_{i-1}(h)| = 1 + a_{i-1}.
$$
This proves the claim. By \eqref{eq:int} we now find $k = c_d + a_d = 2d$ and $b_i = k - a_i - c_i = 2(d-i)$. 
By \cite{BCN}[Corollary 9.2.5], $\Ga$ is the Hamming graph $H(d,3)$. \qed

%%%%%%%%%%%%%%%%%%%%%%%%%%%%%%%%%%%%%%%%%%%%%%%%%%%%%%%%%%%%%%%%%%%%%%%%%%%%%%%%%%%%%%%%%%%%
\section{Case $[G:H] = 3$ and $o(s) \geq 6$}
\label{sec:case2}

In this section we consider the distance-regular graphs $\Ga$, where $\Ga$ is as in Notation~\ref{blank}, for which $[G:H] = 3$ and $o(s) \ge 6$.

\begin{lemma}
\label{le:prop:1}
Let $\G$ be as in Notation~\ref{blank} and suppose $[G:H] = 3$ and $o(s) \ge 6$.
If $\G$ is distance-regular, then $a_1=0$, $c_2=3$ and $o(s)=6$.
\end{lemma}
\proof
Observe that, by Lemma~\ref{lem:easy}, adjacent vertices $1$ and $s$ have no common neighbours, 
implying that $a_1=0$. By Lemma~\ref{lem:c2}, we have that $c_2 \ge 2$, and so since $\partial(1,s^2)=2$, the vertices $1$ and $s^2$
must have at least two common neighbours. One of these common neighbours is $s$
and by Lemma~\ref{lem:easy} the only other two vertices which could be common neighbours
of $1$ and $s^2$ are $s^3$ and $s^{-1}$. However, in both cases we obtain that
$\{s^3, s^{-3}\} \subseteq S$. This shows that $s$, $s^3$ and $s^{-1}$ are all common 
neighbours of $1$ and $s^2$, implying $c_2=3$.

\smallskip \noindent
Consider now the vertices $1$ and $s^4$. Since $o(s) > 3$ and $s^3 \in S$, Lemma~\ref{lem:easy} implies
that $\partial(1,s^4)=2$. Thus $N(1) \cap N(s^4)$ consists of three vertices, two of which are $s$ and $s^3$.
Denote the third by $x$. By Lemma~\ref{lem:easy}, $x \in Hs^2$. But then on one hand $x=s^{-1}$ and on the other hand 
$x=s^5$. It follows that $o(s)=6$. \qed

\begin{proposition}
\label{prop:2}
Let $\G$ be as in Notation~\ref{blank} and suppose $[G:H] = 3$ with $o(s) \ge 6$.
If $\G$ is distance-regular, then $\G$ is isomorphic to $K_{3,3}$.
\end{proposition}
\proof
We first show that $G$ is generated by $s$.
Suppose on the contrary that $G \ne \la s \ra$ and pick $t \in S \setminus \la s \ra$.
Consider the vertices $1$ and $ts$. Since $ts \in Hs$ we have that $ts \notin S$, and so $\partial(1,ts)=2$. By Lemma~\ref{lem:easy}
the only common neighbours of $1$ and $ts$ are $s$ and $t$ (recall that $t \notin \la s \ra$). But this contradicts Lemma~\ref{le:prop:1} 
which states that the intersection number $c_2$ of $\G$ is $3$.
Therefore $G = \la s \ra$, and so $a_1 = 0$ and $c_2=3$ implies $S=\{s,s^3,s^5\}$. Thus $\G$ is the complete bipartite graph $K_{3,3}$. \qed

%%%%%%%%%%%%%%%%%%%%%%%%%%%%%%%%%%%%%%%%%%%%%%%%%%%%%%%%%%%%%%%%%%%%%%%%%%%%%%%%%%%%%%%%%%%%
\section{Case $[G:H] = 2$ and $o(s) \ge 6$}
\label{sec:case3}

Let $\Ga$ be as in Notation~\ref{blank}. Based on the results of the previous sections we are left with the possibility that $H$ is an index two subgroup in $G$. 
In this section we deal with the case when the order of $s$ is at least $6$. The most difficult case, that is when $s$ is of order $4$, is dealt with in the next section. 
We start with the following result, which holds also for the case $o(s) = 4$.

\begin{lemma}
\label{lem:lambda}
Let $\G$ be as in Notation~\ref{blank} and suppose $[G:H] = 2$ with $o(s) \ge 4$.
If $\Ga$ is distance-regular, then the following holds.
\begin{itemize}
\item[{\rm (i)}]
We have $a_1 \in \{0,2\}$. Moreover, $a_1=2$ if and only if $s^2 \in S$. 
\item[{\rm (ii)}]
If $G$ is not generated by $s$, then $c_2 \in \{2,4\}$.
\end{itemize}
\end{lemma}
\proof
(i) Consider adjacent vertices $1$ and $s$ of $\Ga$. By Lemma~\ref{lem:easy}, the only two 
neighbours of $s$ in $H$ are $1$ and $s^2$, and the only two neighbours of $1$ in $Hs$
are $s$ and $s^{-1}$. Therefore, the common neighbours of $1$ and $s$ are contained in 
the set $\{s^2,s^{-1}\}$. Observe that $1$ is adjacent to $s^2$ or $s$ is adjacent to $s^{-1}$
if and only if $s^2 \in S$. In this case $1$ and $s$ are adjacent to both $s^2$ and $s^{-1}$. 
This shows that $1$ and $s$ either have no common neighbours or they have two common neighbours.
This proves (i). 

\smallskip \noindent
(ii) Since $G$ is not generated by $s$, there exists $t \in S \setminus \la s \ra$. 
Observe that $\partial(1,ts)=2$. Let us determine the set $C$ of common neighbours of $1$ and $ts$.
Note that $s,t \in C$ and that, by Lemma~\ref{lem:easy}, the only other possible 
vertices that could be in $C$ are $s^{-1}$ and $ts^2$. Since $s^{-1} \in N(ts)$ if and only if
$ts^2 \in S$, we find that $c_2 \in \{2,4\}$, as claimed. \qed

We first deal with the graphs without triangles, that is with the possibility $a_1 = 0$.

\begin{proposition}
\label{prop:not0}
Let $\G$ be as in Notation~\ref{blank} and suppose $[G:H] = 2$ with $o(s) \ge 6$.
If $\Ga$ is distance-regular with $a_1=0$, then $\G$ is isomorphic to 
$K_{6,6} - 6K_2$, the complete bipartite graph $K_{6,6}$ minus a 1-factor. 
\end{proposition}
\proof
Suppose that $G$ is generated by $s$ so that $G$ is a cyclic group. By \cite{MP}[Theorem 1.2]
(since $k \ge 3$ and $a_1=0$) $\G$ is a complete bipartite graph or a complete 
bipartite graph minus a 1-factor. In particular, $\G$ is bipartite. 
By Lemma~\ref{lem:easy} the only neighbours of $1$ in $Hs$ are $s$ and $s^{-1}$, 
and so $k \ge 3$ implies that $1$ has a neighbour $h$ in $H$.
Since $G=\la s \ra$, we have $h=s^{2\ell}$ for some $\ell$, and so $\G$ has an odd cycle, a contradiction.
We can thus assume that there exists $t \in S \setminus \la s \ra$. By Lemma~\ref{lem:lambda}, $c_2 \in \{2,4\}$.

\medskip \noindent
Consider now the vertices $1$ and $s^2$. Since, by Lemma~\ref{lem:lambda}, $s^2 \notin S$, 
we have that $\partial(1,s^2)=2$. It follows that $1$ and $s^2$ have at least two common neighbours, one of them being $s$. 
Since $o(s) \ge 6$, $s$ is the only common neighbour of $1$ and $s^2$ contained in $Hs$. There thus must exist at least one common neighbor of $1$ and $s^2$, contained in $H$, say $x$. Since $1$ is adjacent with $x$, we have $\{x,x^{-1}\} \subseteq S$ and since $s^2$ is also adjacent with $x$, we have $\{s^{-2} x, s^2 x^{-1}\} \subseteq S$. This shows that also $s^2 x^{-1}$ is a common neighbour of $1$ and $s^2$. Now, if $c_2 = 2$, this implies $x^2=s^2$, but then $s, s^{-1}, x, x^{-1} \in N(1) \cap N(s^{-1}x)$.
Since $\partial(1,s^{-1}x)=2$, this is impossible. By Lemma~\ref{lem:lambda} we thus have that $c_2=4$.

\medskip \noindent
Recall that $t \in S \setminus \la s \ra$. We claim that $t s^{2\ell} \in S$ for $0 \le \ell \le \frac{o(s)-2}{2}$. Note that the claim is true for $\ell = 0$ since $t \in S$. Assume now that $t s^{2\ell} \in S$ for some $0 \le \ell \le \frac{o(s)-4}{2}$. Since $t \notin \la s \ra$, this implies that $\partial(1,ts^{2\ell+1})=2$, and thus the fact that $c_2 = 4$ implies that $N(1) \cap N(ts^{2\ell+1}) = \{s, s^{-1}, ts^{2\ell}, ts^{2\ell+2} \}$. But then $ts^{2\ell+2} \in S$, which, by induction, proves the claim.

\medskip \noindent
We thus have that $\{s, t, ts^2, ts^4, \ldots, ts^{o(s)-2} \} \subseteq N(1) \cap N(s^2)$. Since $c_2=4$
this shows that $o(s)=6$. Moreover, as $t \in S \setminus \la s \ra$ was arbitrary, we must have $S = \{s, s^{-1}, t, ts^2, ts^4\}$, and so $\G$ is clearly isomorphic to $K_{6,6} - 6K_2$.
\qed

We now consider the case $a_1=2$.

\begin{proposition}
\label{prop:a12}
Let $\G$ be as in Notation~\ref{blank} and suppose $[G:H] = 2$ with $o(s) \ge 6$.
If $\Ga$ is distance-regular with $a_1=2$, then $\G$ is isomorphic to $K_{2,2,2}$. 
\end{proposition}
\proof
Recall that, by Lemma~\ref{lem:lambda}, $s^2 \in S$.
Since $1$ and $s^2$ are adjacent, they have two common neighbours, one of them being $s$.
The other one, denote it by $x$, is contained in $H$ by Lemma~\ref{lem:easy}. 
As in the proof of Proposition~\ref{prop:not0} we find that $\{x,x^{-1}, s^{-2}x, s^2 x^{-1}\} \subseteq S$ and that 
consequently $x^2 = s^2$ must hold. Since $x \notin \{1,s^2\}$ Lemma~\ref{lem:easy} implies that $\partial(s,x) = 2$. As $x,x^{-1} \in S$, 
we have that $N(s) \cap N(x) = \{1,s^2,xs,xs^{-1}\}$, and so it follows that $c_2 = 4$. Now, let $t \in S \setminus \{s,s^{-1}, s^2\}$ 
be an arbitrary neighbor of $1$ in $H$, other than $s^2$. Then by Lemma~\ref{lem:easy} we have that $\partial(s,t) = 2$, and so $c_2 = 4$ 
implies that $N(s) \cap N(t) = \{1,s^2,ts,ts^{-1}\}$. In particular, $t$ is a common neighbor of $1$ and $s^2$. But then $t = x$, and so 
$S = \{s,s^{-1},s^2,x\}$. Since $o(s) \geq 6$ it follows that $s^{-2} = x$, and so $x^2 = s^2$ implies $o(s) = 6$. It is now clear that $\G$ is isomorphic 
to the complete multipartite graph $K_{2,2,2}$. 
\qed

%%%%%%%%%%%%%%%%%%%%%%%%%%%%%%%%%%%%%%%%%%%%%%%%%%%%%%%%%%%%%%%%%%%%%%%%%%%%%%%%%%%%%%%%%%%%
\section{Case $[G:H] = 2$ and $o(s) = 4$}
\label{sec:case4}

Let $\Ga$ be as in Notation~\ref{blank}. In this section we consider the case when $[G:H] = 2$ and $o(s) = 4$. 
In the case that $G = \la s \ra$ we have that $\G$ is isomorphic either to $C_4$ or $K_4$. For the rest of this section we will therefore assume that $s$ does
not generate $G$. By Lemma~\ref{lem:lambda}, we have $a_1 \in \{0,2\}$ and $c_2 \in \{2,4\}$.

We first define a certain matrix $P$. The matrix $P$ has rows and columns indexed by the elements of $V_0=V(\G_0)$, having $x,y$ entry 
\begin{equation}
\label{eq:P}
  P_{x, y} = \left\{ \begin{array}{lll}
                 1 & \hbox{if } \; y = x s^2,    & \\
                   &                                   & (x,y \in V_0). \\
                 0 & \hbox{otherwise } & \end{array} \right. 
\end{equation}
In what follows we abbreviate $Px = x s^2$, we let $A$ be the adjacency matrix of $\G_0$, 
and we let $I$ denote the identity matrix of dimension $|V_0| \times |V_0|$.
We denote the valency of $\G_0$ by $k$ (and so the valency of $\G$ is $k+2$)
and the $i$-th distance matrix of $\G$ by $A_i$ (the matrices $A_i$ were defined in Section~\ref{sec:prelim}).
We denote the intersection numbers of $\G$ with $a_i, b_i, c_i$, and the intersection numbers of $\G_0$ (if they exist) by $a_i^0, b_i^0, c_i^0$.
Recall that 
\begin{equation}
\label{bm}
  A_1 A_i = b_{i-1}A_{i-1} + a_i A_i + c_{i+1} A_{i+1} \quad  (0 \le i \le d),
\end{equation}
where the matrices $A_{-1}$ and $A_{d+1}$ are taken to be the zero matrices and $b_{-1} = c_{d+1} = 0$.
For brevity we denote the distance functions of $\G$ and $\G_0$ by  $\partial = \partial_{\G}$ and $\partial_0 = \partial_{\G_0}$. Moreover, the $i$-th sphere in $\G_0$, centred at a vertex $x \in V_0$, will be denoted by $N^0_i(x)$. We abbreviate $N^0(x) = N^0_1(x)$. 

\begin{lemma}
\label{lem41}
Let $\G$ and $\G_0$ be as in Notation~\ref{blank} and suppose $[G:H] = 2$ with $o(s) = 4$. Let $P$ be the matrix defined in \eqref{eq:P}. 
Then the following holds.
\begin{itemize}
\item[{\rm (i)}] $P^2 = I$.
\item[{\rm (ii)}]
Let $M$ be any matrix with rows and columns indexed by the elements of $V_0$ and pick $x,y \in V_0$.
Then $(MP)_{x,y} = M_{x,Py}$ and $(PM)_{x,y} = M_{Px,y}$.
\item[{\rm (iii)}]
$AP =PA$.
\item[{\rm (iv)}]
The mapping which interchanges $x$ and $Px$ for every $x \in V_0$ is an automorphism of the graph $\G_0$.
\end{itemize}
\end{lemma}
\proof
(i) Since $o(s) = 4$, we have that $P(Px) = x$ for any $x \in V_0$. It follows that $P^2 = I$.

\smallskip \noindent
(ii) Since $P^2 = I$ by (i), the definition of $P$ implies 
$$
  (MP)_{x,y} = \sum_{z \in V_0} M_{x,z} P_{z,y} = M_{x,Py}
$$
and
$$
  (PM)_{x,y} = \sum_{z \in V_0} P_{x,z} M_{z,y} = M_{Px,y}.
$$
(iii) 
Pick $x,y \in V_0$. By (ii) above we have $(AP)_{x,y} = A_{x,Py}$ and $(PA)_{x,y} = A_{Px,y}$. The vertices $x$ and $Py = ys^2$ are adjacent if and only if $x^{-1}ys^2 \in S$, which, as $o(s) = 4$, occurs if and only if the vertices $Px = xs^2$ and $y$ are adjacent. The result follows. 

\smallskip \noindent
(iv) Since this is multiplication by $s^2$, this is clear (of course, (iv) also follows from (iii) above). \qed

We remark that, since $Px = xs^2$, the distance $\partial_0(x,Px)$ does not depend on the choice of $x \in V_0$.

%%%% subsection a_1 = 0
\subsection{Subcase $a_1=0$}

Let us first consider the subcase $a_1=0$, that is the subcase when $s^2 \not \in S$.
The distance $\partial_0(x,Px)$ is thus at least $2$ for each $x \in V_0$ (of course $\partial(x,Px) = 2$).  Recall that, by Lemma~\ref{lem:easy}, we have that $x \in H$ and $ys \in Hs$ are adjacent if and only if $y \in \{x, xs^2\} = \{x, Px\}$, and so for an appropriate ordering  
of the vertices of $\G$ (in which all the vertices of $H$ precede all the vertices of $Hs$), the 
adjacency matrix $A_1$ of $\G$ is 
\begin{equation}
\label{adj_mat}
  \left( \begin{array}{cc}
    A   & P+I \cr
    P+I & A \cr
  \end{array}\right).
\end{equation}

\begin{lemma}
\label{lem43}
Let $\G$ be as in Notation~\ref{blank} and suppose $[G:H] = 2$ with $o(s) = 4$.
Assume $\Ga$ is distance-regular with $a_1=0$ and denote $d = \partial_0(x,Px)$ for some
(an thus every) $x \in V(\G)$. Then the following {\rm (i)--(iv)} hold for $1 \le i \le \lfloor d/2 \rfloor$.
\begin{itemize}
\item[{\rm (i)}]
$a_i^0$ exists and $a_i = a_i^0 = 0$.
\item[{\rm (ii)}]
$c_i^0$ exists and $c_i = c_i^0 = i$.
\item[{\rm (iii)}]
There exist polynomials $p_{i-1}(\lambda)$, $q_{i-1}(\lambda)$ and $r_{i-2}(\lambda)$
of degrees $i-1$, $i-1$ and $i-2$, respectively (where the polynomial $r_{-1}(\lambda)$ is taken to be the zero polynomial), such that the $(i+1)$-th distance matrix $A_i$ of $\G$ is given by
\begin{equation}
\label{eq:dm}
  c_{i+1} A_{i+1} = \left( \begin{array}{cc}
    {A^{i+1} + p_{i-1}(A) + q_{i-1}(A)P \over i!} & {((i+1)A^i + r_{i-2}(A))(P+I) \over i!} \cr
    {((i+1)A^i + r_{i-2}(A))(P+I) \over i!} & {A^{i+1} + p_{i-1}(A) + q_{i-1}(A)P \over i!} \cr
  \end{array}\right).
\end{equation}
Moreover, the leading coefficient of $q_{i-1}(\lambda)$ is $i(i+1)$. 
\item[{\rm (iv)}]
For any $x,y \in V_0$ with $\partial_0(x,y) = i$ we have that $\partial(x,y) =  \partial(xs,ys) = i$ and $\partial(x,ys) = \partial(xs,y) = i+1$.
\end{itemize}
\end{lemma}
\proof
Observe that since $a_1 = 0$ we have that $d \geq 2$ (in fact $d \geq 3$ unless $c_2 = 4$).
Fix $x \in V(\G)$. We prove the lemma by induction on $i$. 
Note that the claims (i) and (ii) hold for $i=1$.
Since $a_1 = 0$, \eqref{bm} yields $c_2A_2 = A_1^2 - (k+2)A_0$, and so claim (iii) holds for $i = 1$ by Lemma~\ref{lem41}. That (iv) holds for $i = 1$ is clear since $\partial_0(x,y) = 1$ implies $y \neq Px$. Assume now that the four claims hold for all $j \in \{1,2, \ldots, i\}$, where $i < \lfloor d/2 \rfloor$.

\medskip \noindent
Pick $y \in V_0$ such that $\partial_0(x,y) = i$. Note that $\partial_0(y, Px) \ge d-i > i + 1$. By induction hypothesis it thus follows that the $(y,xs)$-entry of $c_{i+1}A_{i + 1}$ is equal to $(i+1) (A^i)_{y, x}/ i!$. Since $c_j^0 = j$ holds for all $1 \le j \le i$ and $\partial_0(y,x)=i$, we clearly have $(A^i)_{y, x} = i!$ (note that $(A^i)_{y, x}$ is the number of walks in $\G_0$ of length $i$ from $y$ to $x$). Therefore $c_{i+1} = i+1$.

\medskip \noindent
Pick now $y \in V_0$ such that $\partial_0(x,y) = i+1$. Note that $\partial_0(y, Px) \ge d - i - 1 > i$. By induction hypothesis it thus follows that the $(y,x)$-entry of $c_{i+1}A_{i + 1}$ is $(A^{i+1})_{y, x} / i!$. Since $\partial_0(x,y) = i+1$, this entry is nonzero, and so $c_{i+1} = i+1$ implies that $(A^{i+1})_{y,x} = (i+1)!$. It follows that $c_{i+1}^0(x,y) = i+1$. Since $y \in V_0$ was an arbitrary vertex with $\partial_0(x,y) = i+1$ and $\G_0$ is vertex-transitive, this shows that $c_{i+1}^0$ exists and $c_{i+1}^0 = i+1$. Observe that we have also proved that $\partial(x,y) = i+1$ holds.

\medskip \noindent
Let us again pick some $y \in V_0$ with $\partial_0(x,y) = i$. Recall that then $\partial_0(y, Px) > i+1$ holds. 
Moreover, by induction hypothesis $\partial(x,y) = i$ and $\partial(y, xs) = i + 1$ hold. 
Therefore, the $(y, xs)$-entries of the matrices $A_i$, $A_{i+1}$ and $A_{i+2}$ are $0$, $1$ and $0$, respectively. By \eqref{bm}, we have 
\begin{equation}
\label{eq:1}
  c_{i+2} A_{i+2} = A_1 A_{i+1} - b_i A_i - a_{i+1} A_{i+1}.
\end{equation}
Inspecting the $(y, xs)$-entry in this equation we thus find that $a_{i+1} = (A_1 A_{i+1})_{y, xs}$. By Lemma~\ref{lem41} and induction hypothesis the upper-right block of the matrix $A_1 A_{i+1}$ is 
\begin{equation}
\label{eq:2}
  {1 \over (i+1)!} \bigg( ((i+1) A^i + r_{i-2}(A))A(P+I) + (A^{i+1} + p_{i-1}(A) + q_{i-1}(A)P)(P+I) \bigg).
\end{equation}
Therefore, the $(y, xs)$-entry of the matrix $A_1 A_{i+1}$ is equal to
$$
  {i+1+1 \over (i+1)!} (A^{i+1})_{y,x}.
$$
Recall that $\partial(x,y) = i$, and so the $(x,y)$-entry of the matrix $A_{i+1}$ is zero. Thus, since $\partial_0(y,Px) > i+1$, induction hypothesis (claim (iii)) implies that 
the $(x,y)$-entry of the matrix $A^{i+1}$ is zero, implying that $a_{i+1}=0$. Of course, this implies that $a_{i+1}^0$ exists and is equal to $0$. 

\medskip \noindent
The required formula for $c_{i+2} A_{i+2}$ and for the leading coefficient of $q_i(\lambda)$ now follows from \eqref{eq:1},
using usual block-matrix multiplication, induction hypothesis, $a_{i+1}=0$ and Lemma~\ref{lem41}. 

\medskip \noindent
To complete the induction step we only need to prove that for $y \in V_0$ such that $\partial_0(x,y) = i+1$ we have that $\partial(x,ys) = \partial(xs,y) = i+2$ (we already proved that $\partial(x,y) = i+1$ whereas $\partial(xs,ys) = i+1$ follows from the fact that multiplication by $s$ is an automorphism). Since $\partial(x,y) = i+1$ and $a_{i+1} = 0$ the only other possibility is $\partial(x,ys) = i$. But in this case $c_{i+1} > c_{i+1}^0(x,y)$ which contradicts (ii).
 \qed

\begin{lemma}
\label{lem44}
Let $\G$ be as in Notation~\ref{blank} and suppose $[G:H] = 2$ with $o(s) = 4$. Assume $\Ga$ is distance-regular with $a_1=0$ and 
denote $d = \partial_0(x,Px)$ for some (an thus every) $x \in V_0$. 
Then the following {\rm (i), (ii)} hold.
\begin{itemize}
\item[{\rm (i)}]
If $d$ is even, then $a_{d/2+1}^0$ exists and $a_{d/2+1} = a_{d/2+1}^0 = 0$, and $c_{d/2+1}=d+2$. Moreover, for any $x, y \in V_0$ such that $\partial_0(x,y) = d/2 + 1$ 
and $\partial_0(Px,y) = d/2 - 1$ we have $c_{d/2+1}^0(x,y) = d/2+1$.
\item[{\rm (ii)}]
If $d$ is odd, then $c_{(d+1)/2}^0$ exists and $c_{(d+1)/2} = c_{(d+1)/2}^0 = (d+1)/2$, and $a_{(d+1)/2} = (d+3)/2$.
Moreover, for any $x, y \in V_0$ such that $\partial_0(x,y) = (d + 1)/2$ 
and $\partial_0(Px,y) = (d - 1)/2$ we have $a_{(d+1)/2}^0(x,y) = 0$.
\end{itemize}
\end{lemma}
\proof
Fix $x \in V_0$. \\
(i) Pick $y \in V_0$ such that $\partial_0(x,y) = \partial_0(Px,y) = d/2$. By Lemma~\ref{lem43}(iv), we have $\partial(y, xs) = d/2+1$, and so the $(y, xs)$-entry of the matrix $A_{d/2+1}$ equals $1$. By Lemma~\ref{lem43}(iii) we thus get
$$
  c_{d/2 + 1} = {(d/2+1)((A^{d/2})_{y,x} + (A^{d/2})_{y, Px}) \over (d/2)!}.
$$
Since $c_i^0 = i$ for $1 \le i \le d/2$, we have that $(A^{d/2})_{y,x} = (A^{d/2})_{y, Px} = (d/2)!$, and so $c_{d/2+1} = d + 2$.

\smallskip \noindent
Consider now the equation \eqref{eq:1} with $i = d/2$. As in the proof of Lemma~\ref{lem43} we find that $a_{d/2 + 1}$ equals $(A_1A_{d/2+1})_{y,xs}$ 
which in turn is equal to the $(y,x)$-entry of the upper-right block of the matrix $A_1A_{d/2 + 1}$, which is as given in \eqref{eq:2} (where $i = d/2$). 
Therefore, the $(y,xs)$-entry of $A_1 A_{d/2+1}$ is equal to
$$
  {d/2+2 \over (d/2+1)!} \left((A^{d/2+1})_{y,x} + (A^{d/2+1})_{y,Px}\right).
$$
As in the proof of Lemma~\ref{lem43} we find that $(A^{d/2+1})_{y,x} = (A^{d/2+1})_{y,Px} = 0$, and so $a_{d/2+1}=0$. 
Of course, this implies that $a_{d/2 + 1}^0$ exists and is equal to $0$.

\medskip \noindent
Pick now $y \in V_0$ such that $\partial_0(x,y) = d/2+1$ and $\partial_0(Px,y) = d/2-1$. Let us compute the $(x,y)$-entry of $A_{d/2+1}$. 
By Lemma~\ref{lem43}(iii), the $(x,y)$-entry of $A_{d/2+1}$ is 
$$
  {(A^{d/2 + 1})_{x,y} + d/2(d/2+1)(A^{d/2-1})_{Px, y} \over 2(d/2+1)!}.
$$
Since $c_i^0 = i$ for $1 \le i \le d/2$, we have $(A^{d/2-1})_{Px,y} = (d/2-1)!$ and $(A^{d/2 + 1})_{x,y} = c_{d/2+1}^0(x,y)(d/2)!$. 
It follows that $c_{d/2+1}^0(x,y) = d/2+1$. 

\medskip \noindent
(ii) Pick $y \in V_0$ such that $\partial_0(y,x)=(d-1)/2$. Then $\partial_0(y,Px) \ge (d+1)/2$ and, by Lemma~\ref{lem43}(iv), $\partial(y, xs)=(d+1)/2$ holds. The $(y,xs)$-entry of the matrix $A_{(d+1)/2}$ is thus equal to $1$.
By Lemma~\ref{lem43}(iii), the $(y,xs)$-entry of $A_{(d+1)/2}$
is equal to 
$$
  {(d+1)/2 (A^{(d-1)/2})_{y,x} \over c_{(d+1)/2} ((d-1)/2)!}.
$$
Since $c_i^0 = i$ for $1 \le i \le (d-1)/2$, we have 
$(A^{(d-1)/2})_{y,x} = ((d-1)/2)!$, and so $c_{(d+1)/2} = (d+1)/2$.

\medskip \noindent
Pick now $y \in V_0$ such that $\partial_0(y,x)=(d+1)/2$ and note that 
$\partial_0(y,Px) \ge (d-1)/2$.
Let us now compute the $(y,x)$-entry of $A_{(d+1)/2}$.
By Lemma~\ref{lem43}(iii) and since $c_{(d+1)/2} = (d+1)/2$, the $(y,x)$-entry of $A_{(d+1)/2}$
is equal to 
$$
  {(A^{(d+1)/2})_{y,x} \over ((d+1)/2)!}.
$$
As $\partial_0(y,x) = (d+1)/2$, we have that $(A^{(d+1)/2})_{y,x} \neq 0$, and so $(A^{(d+1)/2})_{y,x}  = ((d+1)/2)!$ must hold.
Since $c_i^0 = i$ for $1 \le i \le (d-1)/2$, we have 
$(A^{(d+1)/2})_{y,x} = c_{(d+1)/2}^0(x,y)((d-1)/2)!$, and so $c_{(d+1)/2}^0(x,y) = (d+1)/2$.
Since $x \in V_0$ and $y \in N_{(d+1)/2}^0(x)$ are arbitrary, the intersection number $c_{(d+1)/2}^0$ exists and is equal to $(d+1)/2$.

\medskip \noindent
Pick now $y \in V_0$ such that $\partial_0(y,x)=(d-1)/2$ and $\partial_0(y,Px) = (d+1)/2$. 
By Lemma~\ref{lem43}(iv) we have $\partial(y, xs)=(d+1)/2$, and so by~\eqref{bm}
$$
  (A_1 A_{(d+1)/2})_{y, xs} = a_{(d+1)/2}.
$$
Using Lemma~\ref{lem43}(iii) we find that the $(y, xs)$-entry of $A_1 A_{(d+1)/2}$ is equal to
$$
  {(d+3)/2 \over ((d+1)/2)!} \Big( (A^{(d+1)/2})_{y,x} + (A^{(d+1)/2})_{y,Px} \Big).
$$
Since $c_i^0 = i$ for $1 \le i \le (d+1)/2$ we have $(A^{(d+1)/2})_{y,Px} = ((d+1)/2)!$.
Similarly as above we find that $(A^{(d+1)/2})_{y,x} =  0$ (note that this also follows 
from the fact that $a_i^0 = 0$ for $1 \le i \le (d-1)/2$). It follows that $a_{(d+1)/2} = (d+3)/2$.

\medskip \noindent
Pick now an arbitrary $y \in V_0$ such that $\partial_0(y,x) = (d+1)/2$ and 
$\partial_0(y,Px) = (d-1)/2$. By the argument above the $(y,x)$-entry of $A_{(d+1)/2}$ is nonzero, and so $\partial(y,x) = (d+1)/2$. Since $a_{(d+1)/2}=(d+3)/2$, \eqref{bm} implies that
$$
  (A_1 A_{(d+1)/2})_{y,x} = (d+3)/2.
$$
Using Lemma~\ref{lem43}(iii) (note that the leading coefficient of $q_{(d-3)/2}$ is $(d-1)(d+1)/4$)
we find that the $(y,x)$-entry of $A_1 A_{(d+1)/2}$ is 
$$
  {1 \over ((d+1)/2)!} \bigg( (d+1) (A^{(d-1)/2})_{y,Px} + (A^{(d+3)/2})_{y,x} + (d-1)(d+1)/4 (A^{(d-1)/2})_{y,Px} \bigg).
$$
Since $c_i^0 = i$ for $1 \le i \le (d+1)/2$ and $a_i^0 = 0$ for $1 \le i \le (d-1)/2$,
we find that the $(y,x)$-entry of $A_1 A_{(d+1)/2}$ is
$a_{(d+1)/2}^0(x,y) + (d+3)/2$. Therefore, $a_{(d+1)/2}^0(x,y) = 0$. \qed

\begin{lemma}
\label{lem45}
Let $\G$ be as in Notation~\ref{blank} and suppose $[G:H] = 2$ with $o(s) = 4$.
Assume $\Ga$ is distance-regular with $a_1=0$ and pick $x \in V_0$. Then every vertex of $V_0$
lies on some geodesics between $x$ and $Px$.
\end{lemma}
\proof
Denote $d = \partial_0(x,Px)$. For nonnegative integers $m,n$ we let $D^m_n = N_m^0(x) \cap N_n^0(Px)$.
Observe that if $|n-m| > d$ or $n+m < d$, then $D^m_n = \emptyset$.

\medskip \noindent  
{\bf Case $d$ is even.}
We first show that $N_{d/2}^0(x) = N_{d/2}^0(Px) = D_{d/2}^{d/2}$. To this end let $y \in N_{d/2}^0(x)$. Then $\partial_0(Px,y) \geq d/2$ while, by Lemma~\ref{lem43}(iv), we have that $\partial(x,y) = d/2$ and $\partial(xs,y) = d/2 + 1$. By Lemma~\ref{lem44}(i) we have that $c_{d/2 + 1} = d+2$. Computing the $(xs,y)$-entry of $c_{d/2+1} A_{d/2+1}$ (using Lemma~\ref{lem43}(iii)) we find that
\begin{equation}
\label{eq:lem45}
	d+2 = \frac{(d/2+1)A^{d/2}_{x,y} + (d/2+1)A^{d/2}_{Px,y}}{(d/2)!}.
\end{equation}
Since $c_i^0 = i$ for all $1 \leq i \leq d/2$ we have that $A^{d/2}_{x,y} = (d/2)!$, and so \eqref{eq:lem45} implies $A^{d/2}_{Px,y} = (d/2)!$. It follows that $y \in N_{d/2}^0(Px)$, and so $N_{d/2}^0(x) \subseteq N_{d/2}^0(Px)$. By Lemma~\ref{lem41}(iv) we obtain $N_{d/2}^0(x) = N_{d/2}^0(Px)$ and the claim follows.

\medskip
\noindent
Next we show that for every $0 \leq i \leq d/2$ we have that $N_{d/2 - i}^0(x) = D^{d/2-i}_{d/2+i}$. We show this by induction on $i$. The base case $i = 0$ was settled in the previous paragraph. Suppose then that the claim holds for all $0 \leq j \leq i$, where $i < d/2$, and let us prove it holds for $i+1$ as well. To this end pick a vertex $y \in N_{d/2-i-1}^0(x)$.
Since $\partial_0(x,y) = d/2-i-1$ and since $c_{d/2-i-1}^0 = d/2-i-1$ and $a_{d/2-i-1}^0 = 0$, the vertex $y$ has $k-d/2+i+1$ neighbors
in the set $N_{d/2-i}^0(x) = D^{d/2-i}_{d/2+i}$. Since this also holds for $y \in D^{d/2-i-1}_{d/2+i+1}$, we have that $k-d/2+i+1 > 0$.  Thus every $y \in N_{d/2-i-1}^0(x)$ has at least one neighbor in $D^{d/2-i}_{d/2+i}$, and so for every such $y$ $\partial_0(Px,y) = d/2+i+1$ must hold. It follows that $N_{d/2-i-1}^0(x) = D^{d/2-i-1}_{d/2+i+1}$, which completes the induction step.
Interchanging the role of $x$ and $Px$ we find that for every $0 \leq i \leq d/2$ we also have that 
$N_{d/2 - i}^0(Px) = D^{d/2+i}_{d/2-i}$.

\medskip
\noindent
To complete the $d$ is even case we now only need to show that $D^{d/2+1}_{d/2+1} = \emptyset$.
Pick $z \in D_{d/2-1}^{d/2+1}$. By Lemma~\ref{lem43}(iii) we have that the $(x,z)$-entry of the matrix $A_{d/2+1}$ is nonzero, and so $\partial(x,z)=d/2+1$. We first show that for all $w \in N(z)$ we have that $\partial(x,w) \leq \frac{d}{2}$. This clearly holds for $w \in D^{d/2}_{d/2}$. Moreover, since $\partial(x,Px) = 2$ this also holds for $w \in D_{d/2-2}^{d/2+2}$. Since $a_{d/2 - 1} = 0$, the only remaining neighbors of $z$ are $zs$ and $(Pz)s$. Since $\partial((Px)s,zs) = \partial(Px,z) = d/2-1 = \partial(x,Pz) = \partial(xs,(Pz)s)$ and $x$ is adjacent to both $xs$ and $(Px)s$, these two vertices are also at distance at most $d/2$ from $x$. It thus follows that $b_{d/2 + 1}=0$, and therefore the diameter of $\G$ is $d/2 + 1$. Now, if there exists $y \in D^{d/2+1}_{d/2+1}$, then the $(x,ys)$-entry of the matrix $A_i$ is zero for all $1 \leq i \leq d/2 + 1$ by Lemma~\ref{lem43}(iii), and so $\partial(x,ys) > d/2+1$, which is impossible. Therefore $D^{d/2+1}_{d/2+1} = \emptyset$.

\medskip \noindent  
{\bf Case $d$ is odd.}
We first show that $N_{(d-1)/2}^0(x) = D_{(d+1)/2}^{(d-1)/2}$. To this end let $y \in N_{(d-1)/2}^0(x)$. Then $\partial_0(Px,y) \geq (d+1)/2$ while, by Lemma~\ref{lem43}(iv), we have that $\partial(x,y) = (d-1)/2$ and $\partial(xs,y) = (d+1)/2$. By Lemma~\ref{lem44}(ii) we have that $a_{(d+1)/2} = (d+3)/2$. By \eqref{bm} the $(xs,y)$-entry of $A_1 A_{(d+1)/2}$ equals $(d+3)/2$, and so Lemma~\ref{lem43}(iii) and Lemma~\ref{lem44}(ii) imply that
\begin{equation}
\label{eq:lem45a}
	\frac{d+3}{2} = \frac{(d+3)/2}{((d+1)/2)!} \Big( A^{(d+1)/2}_{x,y} + A^{(d+1)/2}_{Px,y} \Big).
\end{equation}
As $\partial(x,y) = (d-1)/2$ and $a_i^0 = 0$ for all $1 \leq i \leq (d-1)/2$, we have that the $(x,y)$-entry of $A^{(d+1)/2}$ is zero, and so $A^{(d+1)/2}_{Px,y}$ is nonzero (and equals $((d+1)/2)!$). Therefore $\partial_0(Px,y) = (d+1)/2$, and the claim follows. Similarly we show that $N_{(d-1)/2}^0(Px) = D_{(d-1)/2}^{(d+1)/2}$.
As in (i) above we now find that for every $0 \leq i \leq (d-1)/2$ we have that $N_{(d-1)/2 - i}^0(x) = D^{(d-1)/2-i}_{(d+1)/2+i}$ and that $N_{(d-1)/2 - i}^0(Px) = D^{(d+1)/2+i}_{(d-1)/2-i}$.

\medskip
\noindent
We next show that $D_{(d+1)/2}^{(d+1)/2} = \emptyset$. Suppose that there exists $y \in D_{(d+1)/2}^{(d+1)/2}$. 
Since $\partial_0(y,Px)=(d+1)/2$ and $D_{(d-1)/2}^{(d+3)/2} = \emptyset$,
$y$ must be adjacent to at least one vertex $z$ in $D_{(d-1)/2}^{(d+1)/2}$. But then $a_{(d+1)/2}^0(x,z) \ne 0$, contradicting Lemma~\ref{lem44}(ii). Therefore, $D_{(d+1)/2}^{(d+1)/2} = \emptyset$.

\medskip \noindent
To complete the $d$ is odd case we now only need to show that $D^{(d+1)/2}_{(d+3)/2} = D^{(d+3)/2}_{(d+1)/2} = \emptyset$. Pick $z \in D_{(d-3)/2}^{(d+3)/2}$. Since $\partial_0(Px,z) = (d-3)/2$, Lemma~\ref{lem43}(iii) implies that the $(x,z)$-entry of the matrix $A_{(d+1)/2}$ is nonzero, and so $\partial(x,z)=(d+1)/2$. 
Similarly as in the $d$ is even case we show that for all $w \in N(z)$ we have that $\partial(w,x) \leq (d+1)/2$.
It follows that $b_{(d+1)/2}=0$ and therefore the diameter of $\G$ is $(d+1)/2$. Now, if there exists $y \in D^{(d+1)/2}_{(d+3)/2} \cup D^{(d+3)/2}_{(d+1)/2}$, then Lemma~\ref{lem43}(iii) implies that either $\partial(Px,y) \geq (d+3)/2$, or $\partial(x,y) \geq (d+3)/2$, a contradiction. Therefore $D^{(d+1)/2}_{(d+3)/2} = D^{(d+3)/2}_{(d+1)/2} = \emptyset$, which completes the proof. \qed

\begin{proposition}
\label{pro:order4a10}
Let $\G$ be as in Notation~\ref{blank} and suppose $[G:H] = 2$ with $o(s) = 4$. 
Pick $x \in V_0$ and denote $d=\partial_0(x,Px)$.
If $\Ga$ is distance-regular with $a_1=0$, then the following (i), (ii) hold.
\begin{itemize}
\item[{\rm (i)}]
$\G_0$ is isomorphic to the $d$-dimensional hypercube $H(d,2)$.
\item[{\rm (ii)}]
$\G$ is isomorphic to the antipodal quotient of 
the $(d+2)$-dimensional hypercube $H(d+2,2)$.
\end{itemize}
\end{proposition}
\proof
(i) Let $k$ denote the valency of $\G_0$. We first show that $k=d$. Indeed, if $d$ is even, pick $y \in D_{d/2}^{d/2}$.
By Lemma~\ref{lem43} we have that $a_{d/2}^0 = 0$ and $c_{d/2}^0=d/2$, and so $y$ has $d/2$ neighbours in $D_{d/2+1}^{d/2-1}$, $d/2$ neighbours in $D_{d/2-1}^{d/2+1}$, and no other neighbours in $\G_0$. This shows that $k=d$. If $d$ is odd, then pick $y \in D_{(d-1)/2}^{(d+1)/2}$. By Lemma~\ref{lem43} and Lemma~\ref{lem44} we have that $a_{(d-1)/2}^0 = 0$ and $c_{(d-1)/2}^0=(d-1)/2$, $c_{(d+1)/2}^0=(d+1)/2$, and so $y$ has $(d-1)/2$ neighbours in $D_{(d-3)/2}^{(d+3)/2}$, $(d+1)/2$ neighbours in $D_{(d+1)/2}^{(d-1)/2}$, and no other neighbours in $\G_0$. This shows that $k=d$.

\medskip \noindent
Combining this with Lemma~\ref{lem43}, Lemma~\ref{lem44} and Lemma~\ref{lem45}, we find that for $y \in N_i^0(x)$ with $0 \le i \le d$ we have $c^0_i(x,y) = i$ and $b^0_i(x,y)=d-i$. Since $x$ is an arbitrary vertex of $\G_0$, it follows that the intersection numbers of $\G_0$ exist and so the intersection array of $\G_0$ is $\{d,d-1,d-2, \ldots, 1 ; 1,2,3, \ldots, d\}$. By \cite[Corollary 9.2.5]{BCN}, $\G_0$ is isomorphic to the $d$-dimensional hypercube.

\medskip \noindent
(ii) Since the valency of $\G_0$ is $d$, the valency of $\G$ is of course $d+2$.
It follows from Lemma~\ref{lem45}, that the diameter of $\G$ is $\lfloor d/2 \rfloor + 1$. 
Combining together Lemma~\ref{lem43} and Lemma~\ref{lem44} we find that the intersection array of $\G$ is 
$$
  \{d+2, d+1, \ldots, d/2+2; 1,2, \ldots, d/2, d+2\}
$$
if $d$ is even, and
$$
  \{d+2, d+1, \ldots, (d+5)/2; 1,2, \ldots, (d-1)/2, (d+1)/2\}
$$
if $d$ is odd.
If $d+2 \ne 6$, then, by \cite[Theorem 9.2.7]{BCN}, $\G$ is isomorphic to the antipodal quotient of the 
$(d+2)$-dimensional hypercube. As for $d+2=6$, \cite[Theorem 9.2.7]{BCN} implies that there are precisely three
nonisomorphic graphs with the intersection array $\{6,5,4;1,2,6\}$, one of them being the antipodal quotient of
the $6$-dimensional hypercube. It is easy to see,  since $\G_0$ is the $4$-dimensional hypercube, that the other two graphs cannot arise in this way.  \qed

\subsection{Subcase $a_1=2$}

Let us now consider the subcase $a_1=2$, that is the subcase when $s^2 \in S$.
Let $\G_1=\cay(H,S \setminus \{s,s^2,s^3\})$, let $V_1 = V(\G_1) = H$ and let $\partial_1$ be the distance function of the graph $\G_1$. We denote the intersection numbers of $\G_1$ (if they exist) by $a_i^1$, $b_i^1$ and $c_i^1$. Note that for any $x \in V(\G_1)$ we have that $\partial_1(x,Px) \geq 2$. Note also that $\G_1$ may be disconnected. For an appropriate ordering of the vertices of $\G$, the adjacency matrix of $\G$ is 
\begin{equation}
\label{adj_mat1}
  \left( \begin{array}{cc}
    A+P   & P+I \cr
    P+I & A+P \cr
  \end{array}\right),
\end{equation}
where $A$ is the adjacency matrix of $\G_1$.

\begin{lemma}
\label{lem51}
Let $\G$ be as in Notation~\ref{blank} and suppose $[G:H] = 2$ with $o(s) = 4$.
Assume $\Ga$ is distance-regular with $a_1=2$. If $\Ga$ is not $K_4$, then $c_2=2$.
\end{lemma}
\proof
Assume that $\G$ is not $K_4$ so that there exists $t \in S \setminus \{s,s^2,s^3\}$. Of course $\partial(1,ts) = 2$. By Lemma~\ref{lem:lambda} we have that $c_2 \in \{2,4\}$. If $c_2 = 4$, then $N(1) \cap N(ts) = \{s,s^3,t,ts^2\}$, implying that $ts^2 \in S$. But then $\{s,s^3,t,ts^2\} \subseteq N(1) \cap N(s^2)$, contradicting $a_1 = 2$.
\qed 

\begin{lemma}
\label{lem52}
Let $\G$ be as in Notation~\ref{blank} and suppose $[G:H] = 2$ with $o(s) = 4$.
If for some (and thus every) $x \in V_1$ the vertices $x$ and $Px$ are not in the same connected component of $\G_1$, then let $d$ be the diameter of one (and thus each) connected component of $\G_1$. Otherwise let $d = \lfloor((\partial_1(x,Px) - 1)/2)\rfloor$.
Assume $\Ga$ is distance-regular with $a_1=2$.
Then the following {\rm (i)--(iv)} hold for $1 \le i \le d$.
\begin{itemize}
\item[{\rm (i)}]
$a_i^1$ exists and $a_i = a_i^1 = 2i$ holds.
\item[{\rm (ii)}]
$c_i^1$ exists and $c_{i} = c_i^1 = i$.
\item[{\rm (iii)}]
There exist polynomials $p_i(\lambda)$, $q_i(\lambda)$ and $r_{i-1}(\lambda)$ 
with degrees $i$, $i$ and $i-1$ respectively, such that the $(i+1)$-th distance matrix $A_i$ of $\G$ is given by
\begin{equation}
\label{eq:dm1}
  c_{i+1} A_{i+1} = \left( \begin{array}{cc}
    {A^{i+1} + p_i(A) + q_i(A)P \over i!} & {((i+1)A^i + r_{i-1}(A))(P+I) \over i!} \cr
    {((i+1)A^i + r_{i-1}(A))(P+I) \over i!} & {A^{i+1} + p_i(A) + q_i(A)P \over i!} \cr
  \end{array}\right).
\end{equation}
Moreover, the leading coefficient of $p_i(\lambda)$, $q_i(\lambda)$ and $r_{i-1}(\lambda)$ are $-(i+1)i$,
$(i+1)$ and $-(i+1)i(i-1)$, respectively.
\item[{\rm (iv)}]
Pick $z,w \in V_1$ with $\partial_1(z,w) = i$. 
Then $\partial(z,w) = \partial(zs,ws) = i$ and $\partial(z,ws) = \partial(zs,w) = i+1$.
\end{itemize}
\end{lemma}
\proof
The proof is by induction on $i$. Using~\eqref{bm} it can be verified that the proposition holds for $i=1$. Assume now that the proposition holds for some $1 \le i \le d-1$ and let us prove it holds also for $i+1$. Pick $x,y \in V_1$ such that $\partial_1(x,y) = i$. By induction hypothesis $\partial(y, xs)=i+1$, and so the $(y,xs)$-entry of $c_{i+1}A_{i+1}$ equals $c_{i+1}$. Hence induction hypothesis implies that $c_{i+1} = (i+1)/(i!)(A^i)_{x,y}$. Since $c^1_j = j$ for all $1 \leq j \leq i$ we thus get $c_{i+1} = i+1$. By~\eqref{bm} we also have that $(A_1A_{i+1})_{y,xs} = a_{i+1}$. Using the fact that $c_{i+1} = i+1$ and induction hypothesis we find that $(A_1A_{i+1})_{y,xs}$ equals the $(y,x)$-entry of the matrix
$$
\frac{1}{(i+1)!}\left( (i+2)A^{i+1} + p_i(A) + q_i(A) + A r_{i-1}(A) + (i+1)A^i + r_{i-1}(A) \right)(P+I).
$$
Recall that since $\partial_1(x,y) = i < d$, we have that $\partial_1(y,Px)>d>i$. Using the induction hypothesis information about the leading coefficients of the polynomials $p_i(\lambda)$, $q_i(\lambda)$ and $r_{i-1}(\lambda)$ we find that $a_{i+1}$ equals the $(y,x)$-entry of the matrix
$$
\frac{1}{(i+1)!}\left( (i+2)A^{i+1} + (i+1)(2-i^2)A^i \right).
$$
By induction hypothesis we have that $c_j^1 = j$ and $a_j^1 = 2j$ for all $1 \leq j \leq i$, and so $(A^i)_{y,x} = i!$ and  $(A^{i+1})_{y,x} = (a^1_1 + a^1_2 + \cdots + a^1_i)i! = i(i+1)!$. It thus follows that $a_{i+1} = 2(i+1)$, as claimed.

\medskip \noindent
By~\eqref{bm} it now follows that $c_{i+2}A_{i+2} = A_1A_{i+1} - b_iA_i - 2(i+1)A_{i+1}$. Using the induction hypothesis information about the leading coefficients of the polynomials $p_i(\lambda)$, $q_i(\lambda)$ and $r_{i-1}(\lambda)$ it can be verified that $c_{i+2}A_{i+2}$ is as in (iii) of the proposition and that the leading coefficients of the polynomials $p_{i+1}(\lambda)$, $q_{i+1}(\lambda)$ and $r_i(\lambda)$ are $-(i+1)(i+2)$, $(i+2)$ and $-i(i+1)(i+2)$, respectively.

\medskip \noindent
Pick now $y \in \G_1$ such that $\partial_1(x,y) = i+1$ and note that $\partial_1(Px,y) > d \geq i+1$. We now compute the $(x,y)$-entry of the matrix $A_{i+1}$. Since $c_{i+1} = i+1$, the induction hypothesis implies it is equal to  $(1/(i+1)!)(A^{i+1})_{x,y} \neq 0$. It follows that $c_{i+1}^1(x,y) = i+1$. Since $x, y \in V_1$ were arbitrary subject to $\partial_1(x,y) = i+1$, the intersection number $c_{i+1}^1$ exists and equals $i+1$. The above argument implies that the $(x,ys)$-entry of the matrix $c_{i+2}A_{i+2}$ equals $(1/(i+1)!)(i+2)(A^{i+1})_{x,y} \neq 0$, implying that $\partial(x,ys) = i+2$, while the $(x,y)$-entry of this matrix, which of course is zero, equals $1/(i+1)!((A^{i+2})_{x,y} + (p_{i+1}(A))_{x,y})$. Hence $0 = (i+1)!(a^1_1 + a^1_2 + \cdots + a^1_i + a^1_{i+1}(x,y)) -(i+1)(i+2)(i+1)!$. Since $a^1_j = 2j$ for all $1 \leq j \leq i$, we thus have that $a_{i+1}^1(x,y) = 2(i+1)$. Since $x, y \in V_1$ were arbitrary subject to $\partial_1(x,y) = i+1$, the intersection number $a_{i+1}^1$ exists and equals $2(i+1)$.

\medskip \noindent
Part (iv) of the lemma is now easily obtained by computing the $(z,w)$-entry and the $(zs,w)$-entry of appropriate distance matrices $A_{i+1}$ and $A_{i+2}$.
\qed

\begin{lemma}
\label{lem53}
Let $\G$ be as in Notation~\ref{blank} and suppose $[G:H] = 2$ with $o(s) = 4$.
Assume $\Ga$ is distance-regular with $a_1=2$. If for some (and thus every) $x \in V_1$ the vertices $x$ and $Px$ are not in the same connected component of $\G_1$, then let $d$ be the diameter of one (and thus each) connected component of $\G_1$. Otherwise let $d = \lfloor((\partial_1(x,Px) - 1)/2)\rfloor$. Then $c_{d+1} = d+1$. Moreover, if $x$ and $Px$ are not in the same component of $\G_1$ then $a_{d+1} = 2(d+1)$.
\end{lemma}
\proof
Pick $x,y \in V_1$ such that $\partial_1(Px,y) > \partial_1(x,y) = d$.
By Lemma~\ref{lem52}(iv) we have that $\partial(y,xs)=d+1$ and 
therefore Lemma~\ref{lem52}(iii) implies that 
$$
  c_{d+1} = \frac{(d+1)(A^d)_{y,x}}{d!}.
$$
Since $c_j^1=j$ for $1 \le j \le d$, we have that $c_{d+1}=d+1$. 
By~\eqref{bm} we thus find that $a_{d+1}$ equals the $(x,ys)$-entry of the matrix $A_1A_{d+1}$. Thus, if $x$ and $Px$ are not in the same component of $\G_1$, then, using Lemma~\ref{lem52}(iii), one can verify that this entry equals $2(d+1)$.
\qed

\begin{lemma}
\label{lem54}
Let $\G$ be as in Notation~\ref{blank} and suppose $[G:H] = 2$ with $o(s) = 4$.
Assume $\Ga$ is distance-regular with $a_1=2$ and pick $x \in V_1$. 
Then $x$ and $Px$ are not in the same connected component of $\G_1$. 
\end{lemma}
\proof
Suppose on the contrary that $x$ and $Px$ are in the same connected component of $\G_1$.
Let $x=x_0, x_1, x_2, \ldots, x_d = Px$ be a shortest path between $x$ and $Px$ in $\G_1$.
Assume first that $d$ is even and consider $y=x_{d/2 + 1}$. By Lemma~\ref{lem52}(iii) the $(x,y)$-entry of the matrix $A_{d/2}$ is nonzero, and so $\partial(x,y)=d/2$. Hence $|N_{d/2}(x) \cap N(y)| = a_{d/2} = d$. But, again by Lemma~\ref{lem52}(iii),
$$
   (N_{d/2-1}^1(Px) \cap N^1(y)) \cup \{ys, (Py)s, x_{d/2}\} \subseteq N_{d/2}(x) \cap N(y)
$$
holds, and so 
$$
   d = a_{d/2} = |N_{d/2}(x) \cap N(y)| \ge 3 + |N_{d/2-1}^1(Px) \cap N^1(y)| = 3 + a_{d/2-1}^1 = d+1,
$$
a contradiction.

\medskip \noindent 
Assume next that $d$ is odd and consider $y = x_{(d+1)/2}$. Using Lemma~\ref{lem52}(iii) we find that $\partial(x,y)=(d+1)/2$, and so Lemma~\ref{lem53} implies $|N_{(d-1)/2}(x) \cap N(y)| = c_{(d+1)/2} = (d+1)/2$. But by Lemma~\ref{lem52} we have that
$$
   (N^1(y) \cap N^1(x_{(d-3)/2})) \cup (N^1(y) \cap N^1_{(d-3)/2}(Px)) \subseteq N(y) \cap N_{(d-1)/2}(x).
$$
As the sets $N^1(y) \cap N^1(x_{(d-3)/2})$ and $N^1(y) \cap N^1_{(d-3)/2}(Px)$ are clearly disjoint, we thus have
$$
  \frac{d+1}{2} = c_{(d+1)/2} = |N_{(d-1)/2}(x) \cap N(y)| \ge $$
$$
|N^1(y) \cap N^1(x_{(d-3)/2})| + |N^1(y) \cap N^1_{(d-3)/2}(Px)| = c_2^1 + c_{(d-1)/2}^1 = \frac{d+3}{2},
$$
a contradiction.
\qed

\begin{proposition}
\label{prop:order4a12}
Let $\G$ be as in Notation~\ref{blank} and suppose $[G:H] = 2$ with $o(s) = 4$.
If $\Ga$ is distance-regular with $a_1=2$, then $\G$ is isomorphic either to the Hamming graph $H(d+1,4)$
or to the Doobs graph $D(n,m), \; n \ge 1$, where $n+2m = d+1$ and where $d+1$ is the diameter of $\G$.
\end{proposition}
\proof
Pick $x \in V_1$. By Lemma~\ref{lem54}, $x$ and $Px$ are not in the same component of $\G_1$.
Since the mapping sanding $z$ to $Pz$ is an automorphism of $\G_1$, there are only two (isomorphic) components 
of $\G_1$. Let $d$ be the diameter of these components. Note that the diameter of $\G$ is at most
$d+1$. By Lemma~\ref{lem53} we have that $c_{d+1} = d+1 \ne 0$, and so the diameter of $\G$ is $d+1$. By Lemma~\ref{lem53} we have that $a_{d+1} = 2(d+1)$ (since $y$ and $Px$ are in different components of $\G_1$). Therefore, by Lemma~\ref{lem52} the intersection numbers of $\G$ are given by $c_i = i$ and $a_i = 2i$ for $1 \le i \le d+1$. This implies that
$b_0 = a_{d+1} + c_{d+1} = 2(d+1) + (d+1) = 3(d+1)$, and therefore 
$b_i = b_0 - a_i - c_i = 3(d+1-i)$. By \cite{BCN}[Corollary 9.2.5], $\G$ is either isomorphic to the Hamming graph $H(d+1,4)$ or to the Doobs graph $D(n,m)$, where $n+2m = d+1$.
Since $\G_1$ consists of two (isomorphic) connected components, $\G$ is clearly isomorphic to the cartesian product of this connected component with the complete graph $K_4$. Since the clique number of the Shrikhande graph is $3$, this shows that in the case when $\G$ is the Doobs graph $D(n,m)$, we have $n \ge 1$. \qed

%%%%%%%%%%%%%%%%%%%%%%%%%%%%%%%%%%%%%%%%%%%%%%%%%%%%%%%%%%%%%%%%%%%%%%%%%%%%%%%

%%%%%%%%%%%%%%%%%%%%%%%%%%%%%%%%%%%%%%%%%%%%%%%%%%%%%%%%%%%%%%%%


\begin{thebibliography}{99}

\begin{footnotesize}

\bibitem{BM}       W.\ G.\ Bridges and R.\ A.\ Mena, Rational circulants with 
                   rational spectra and cyclic strongly regular graphs, 
                   {\em Ars Combin.}, {\bf 8} (1979), 143--161.
\bibitem{BCN}      A.\ E.\ Brouwer, A.\ M.\ Cohen and A.\ Neumaier, 
                   {\em Distance-regular graphs},
                   Springer-Verlag, New York (1998).
\bibitem{MaSurvey} S.\ L.\ Ma, A Survey of Partial Difference Sets,
                   {\em Designs, Codes and Cryptography}, {\bf 4} (1994), 221--261.
\bibitem{Ma1}      S.\ L.\ Ma, Partial Difference Sets,
                   {\em Discrete math.}, {\bf 52} (1984), 75--89.
\bibitem{DM}       D.\ Maru\v si\v c, Strong regularity and circulant graphs,
                   {\em Discrete math.}, {\bf 78} (1989), 119--125.
\bibitem{MP}       \v{S}.~Miklavi\v{c} and P. Poto\v{c}nik, Distance-regular circulants, 
                   {\em Europ. J. Combin.} {\bf 24} (2003), 777--784.
\bibitem{MP1}      \v{S}.\ Miklavi\v{c} and P.\ Poto\v cnik,
                   Distance-regular Cayley graphs on dihedral gro\-ups,
                   {\it J. Combin. Theory Ser.~B}, {\bf 97} (2007), 14--33.
\bibitem{Mu}       M. Muzychuk, Strongly regular Cayley graphs over the group $\ZZ_{p^n} \oplus \ZZ_{p^n}$,
                   {\em Discrete Math.}, {\bf 305} (2005), 219--239.
\bibitem{SS59} S.~S.~Shrikhande, The uniqueness of the $L_2$ association scheme, {\em Ann. of Math. Stat.}, {\bf  30} (1959), 781--798.
\end{footnotesize}
\end{thebibliography}
\end{document}